\documentclass[12pt]{article}

\usepackage[T2A]{fontenc}
\usepackage[utf8]{inputenc}

\usepackage{amsfonts}
\usepackage{amssymb}
\usepackage{amsmath}
\usepackage{amsthm}
\usepackage{mathrsfs}

\usepackage{tikz}
\usetikzlibrary{calc,arrows,decorations.pathreplacing,fadings,3d,positioning}

\def \le {\leqslant}
\def \ge {\geqslant}

\theoremstyle{plain}
 
\begin{document}
\begin{Large}
\centerline{Weak approximations, Diophantine exponents}
\centerline{ and two-dimensional lattices}
\end{Large}
\vskip+0.5cm

\centerline{Nikolay Moshchevitin\footnote{Technische Universität Wien, research is supported by Austrian Science Fund (FWF), Forschungsprojekt PAT1961524.}}

\vskip+1cm
 
 In the present paper we study the properties of Diophantine exponents of lattices and so-called related "weak" uniform approximations introduced in recent papers \cite{G2,G1}, in the simplest two-dimensional case.  In contrast to the multidimensional case, in the two-dimensional case we can use a powerful tool of continued fractions. We develop an analog of Jarn\'{\i}k's theory  (see original paper \cite{J} as well as a general result \cite{Mamo}
 and an extended survey \cite{G0}) dealing with inequalities between the ordinary and uniform Diophantine exponents, which turned out to be related to mutual behaviour of irrationality measure functions for two real numbers
 (see  the main result from \cite{KM}  and its extensions \cite{D, M1, SRu, S}).
 
 \vskip+0.3cm
 The structure of our paper is as follows. In Section 1 we recall the basic  properties of the ordinary irrationality measure function  for one number, discuss the connection with continued fractions,  introduce the concept of  weak approximations and  mutual behaviour of two irrationality measure functions. Moreover, here we recall the definitions of Diophantine exponents for approximations to one  number  and introduce exponents for the minima of  two irrationality measure functions. Finally, we explain the connection with Diophantine exponents for two-dimensional lattices. In Section 2 we formulate our main results dealing with inequalities between different exponents for one and two numbers.
 The most important result here is Theorem 3 which gives a lower bound for the minimum of two ordinary Diophantione exponents  $\omega(\theta), \omega(\eta)$ for two real numbers $\theta,\eta$ in terms of mutual uniform Diophantine exponents and its corollary Theorem 4, which gives a lower bound for the ordinary Diophantine exponent $\omega(\Lambda)$ of two-dimensional Lattice $\Lambda$ in terms of its weak uniform exponent  $\overline{\omega}(\Lambda)$.
 In Section 3 we formulate and prove    a lemma about mutual behaviour of two piecewise constant functions. This lemma is not related to the Diophantine nature of the problem under consideration and does not use continued fractions. In Section 4 we finalise the proofs of our main theorems. The last Section 5 is devoted to explanation of optimality of the inequalities obtained in our theorems.
 
\vskip+0.3cm
{\bf 1. Irrationality measure functions}
\vskip+0.3cm
Let $\theta$ be an irrational real number. We deal with its representation as a continued fraction
\begin{equation}\label{1m}
\theta=
   [a_0;a_1,a_2,...,a_\nu,...]
   =   a_0  +
\frac{1}{\displaystyle{a_1+\frac{1}{\displaystyle{a_2 + \cdots+
\frac{1}{\displaystyle{a_\nu + ...
{} }}}}}}
   , a_0\in \mathbb{Z}, a_\nu \in \mathbb{Z}_+, \nu=1,2,3,... 
\end{equation}
Let
\begin{equation}\label{1mx}
\frak{Q}:\,q_1<q_2<...<q_\nu<q_{\nu+1}<...
\end{equation}
be the sequences of convergents' denominators for $\theta$.
It is well known that 
\begin{equation}\label{2m}
\frac{1}{2q_{\nu+1}}<
 ||q_\nu \theta||  < \frac{1}{q_{\nu+1}};\,\,\,\,\,
\frac{1}{(a_{\nu+1}+2)}<
q_\nu \cdot ||q_\nu \theta||  < \frac{1}{a_{\nu+1}};\,\,\,\,\,
q_\nu \cdot ||q_\nu \theta|| \asymp \frac{q_\nu}{q_{\nu+1}}.
\end{equation}
Here $||\cdot ||$ stands for the distance to the nearest integer, and everywhere in the paper we use Vinorgadov symbols $\asymp,\gg,\ll$.
 Basic facts about continued fractions and in particular proofs of inequalities (\ref{2m}) one can find in   \cite{KI}.
\vskip+0.3cm

{\bf 1.1. Ordinary rational approximations}
\vskip+0.3cm
For $\theta \in \mathbb{R}$ 
we consider the ordinary irrationality measure function defined as
\begin{equation}\label{3m}
\psi_\theta (t) = \min_{q\in \mathbb{Z}_+:q\le t} ||q\theta||.
\end{equation}
We should note that 
\begin{equation}\label{4m}
\psi_\theta (t)  \le \frac{1}{t} \,\,\,\,\,\text{for any}\,\,\,\,\, t\ge 1.
\end{equation}
By Lagrange's theorem we have an explicit formula
$$
\psi_\theta (t) = ||q_\nu \theta|| \,\,\,\,\,
\text{for all}\,\,\, t\,\,\,\text{from  the interval}
\,\,\,\,\,
q_\nu \le t < q_{\nu+1},
$$
so $\psi_\theta (t) $ is a  decreasing piecewise constant function which is non-continuous in every point
$ t= q_\nu, \nu = 1,2,3,...$\, .

The
{\it ordinary Diophantine exponent} $\omega (\theta)$  is defined as 
$$
\omega (\theta) = \sup \{ \gamma\in \mathbb{R}:\,\,
\liminf_{t\to \infty}
t^\gamma \cdot \psi_\theta (t) <\infty\}.
$$
In other words $\omega (\theta)$ is the supremum of those $\gamma$ for which the inequality
$$
||q\theta || < q^{-\gamma} 
$$
has infinitely many solutions in positive integers $q$.
From (\ref{4m}) it follows that $ \omega(\theta) \ge 1$. Moreover, 
given  $\omega \in [1, +\infty]$,  it is easy to construct by means of continued fraction representation and   (\ref{2m}) 
a number
$\theta \in \mathbb{R}\setminus \mathbb{Q}$ such that $ \omega(\theta) = \omega$.

\vskip+0.3cm
{\bf 1.2. Weak approximations}
\vskip+0.3cm
We define 
{\it weak irrationality measure function} as
$$
\upsilon_\theta (t) = \min_{q\in \mathbb{Z}_+:q\le t} q\cdot ||q\theta||.
$$
It is clear that 
\begin{equation}\label{ppq}
\upsilon_\theta (t) \le t \cdot \psi_\theta (t).
\end{equation}
 For  $
\upsilon_\theta (t)$ we have an alternative formula 
\begin{equation}\label{alt}
\upsilon_\theta (t) = \min_{\nu: \,\, q_\nu \le t } \,\, q_\nu\cdot ||q_\nu\theta||.
\end{equation}
Of course, $
\upsilon_\theta (t) $ is a constant function on the interval $q_\nu \le t < q_{\nu+1}$,
but it can happen that $
\upsilon_\theta (t)$ is continuous in a point of the form $ t = q_\nu$ for some values of $\nu$. 
However, if $
\upsilon_\theta (t) $ is not continuous at point $ t= q_\nu$, then
\begin{equation}\label{alto}
\upsilon_\theta (t) = q_\nu\cdot ||q_\nu\theta||\,\,\,\,\,
\text{for all}\,\,\, t\,\,\,\text{from  the interval}
\,\,\,\,\,
q_\nu \le t < q_{\nu+1}.
\end{equation}

%Let us consider only those elements of the  sequence $\frak{Q}$ where $
%\upsilon_\theta (t)$ 
%is not continuous. These elements form a subsequence $\frak{Q}'\subset \frak{Q}$. Let this subsequence be
%$$
%\frak{Q}':\,q_1'<q_2'<...<q_k'<q_{k+1}'<... \, ,
%$$
%so $q_k' = q_{\nu_k}$ for some increasing sequence of indices $ \nu_k \ge k$.

Following  \cite{G1}\footnote{Although we are not sure that the definition  of weak Diophantine exponent is given in the most convenient way, 
in our paper
we decided not to introduce new objects and to keep the original notation from \cite{G1} with 
$\gamma-1$ in exponent in (\ref{food}).
}
we define 
the {\it weak uniform Diophantine exponent} $\overline{\omega}(\theta)$ as
\begin{equation}\label{food}
\overline{\omega} (\theta) =\sup\{\gamma \in \mathbb{R}:
\,\,\,
\limsup_{t\to \infty} t^{\gamma-1} \upsilon_\theta (t)<\infty
\}.
\end{equation}
In other words $\omega (\theta)$ is the supremum of those $\gamma$ for which the  system of inequalities
$$
q\cdot ||q\theta || < t^{1-\gamma} ,\,\,\,\,\, 1\le q \le t
$$
is  solvable in integers for all $t$  large enough. It was shown in \cite{G1}  that 
for any irrational $\theta$ one has $ 1\le \overline{\omega} (\theta)\le 2$ and that  for any $\overline{\omega} \in [1,2]$ there exists $\theta $ such that 
$\overline{\omega} (\theta)= \overline{\omega} $. As it was mentioned in \cite{G1}, obviously we have inequality
\begin{equation}\label{j}
 {\omega} (\theta)\ge \overline{\omega} (\theta).
\end{equation}

{\bf 1.3 Two irrationality measure functions}

\vskip+0.3cm

For two irrational numbers $\theta$ and $\eta$ we consider functions
$$
\psi_{\theta,\eta}  (t)= \min ( \psi_\theta (t), \psi_\eta (t))
$$
and
$$
\upsilon_{\theta,\eta}  (t)= \min ( \upsilon_\theta (t), \upsilon_\eta (t)).
$$
The analysis of mutual behaviour of the  functions $\psi_\theta (t)$ and $\psi_\eta (t)$ was probably initiated in \cite{KM}.
There are several recent papers devotes to this topic, see for example \cite{D,M1,S,SRu}. 
In \cite{Kl} the authors developing a classical construction by Khintchine \cite{HIN} proved a general result which in particular shows  that for any function $\phi(t)$ decreasing to zero there exist two numbers $\theta$ and $\eta$ linearly independent  together with $1$ over $\mathbb{Q}$ such that  for both functions we have
$$
\psi_{\theta,\eta}  (t)\le \phi(t) \,\,\,\,\,\text{and}
\,\,\,\,\,
\upsilon_{\theta,\eta}  (t) \le \phi(t) 
$$
for all $t$ large enough. Such a result for $
\psi_{\theta,\eta}  (t)$ was used in \cite{M0} in a hidden form  (see construction  from the proof of  Theorem 3  from \cite{M0}). 
A similar result for $ \upsilon_{\theta,\eta}  (t)$ was observed in \cite{G2}.

 In the present paper we are interested in  the uniform Diophantine exponents related to the behaviour of functions
 $
 \psi_{\theta,\eta}  (t)
 $
 and
 $
\upsilon_{\theta,\eta}  (t)$. So we define
$$
\varpi_\psi (\theta,\eta) 
=\sup\{\gamma \in \mathbb{R}:
\,\,\,
\limsup_{t\to \infty} t^{\gamma} \psi_{\theta,\eta}  (t) <\infty
\}
$$
and
$$
 {\varpi}_\upsilon (\theta,\eta) 
=\sup\{\gamma \in \mathbb{R}:
\,\,\,
\limsup_{t\to \infty} t^{\gamma-1} \upsilon_{\theta,\eta}  (t) <\infty
\}.
 $$
 Both exponents $\varpi_\psi(\theta,\eta) $ and ${\varpi}_\upsilon (\theta,\eta) $  can attain any values in the interval
 $[1,+\infty]$.
 From (\ref{ppq}) we see that 
\begin{equation}\label{ppq1}
1\le 
\varpi_\psi (\theta,\eta) 
\le
{\varpi}_\upsilon (\theta,\eta) .
\end{equation}

\vskip+0.3cm

{\bf 1.4. Two-dimensional lattices}
\vskip+0.3cm

We consider two-dimensional lattice
$$
\Lambda =
\Lambda_A = A\mathbb{Z}^2,\,\,\,\,\,
A=
\left(
\begin{array}{cc}
a_{1,1}& a_{1,2}\cr
a_{2,1}&a_{2,2}
\end{array}
\right),\,\,\,\,\,
{\rm det}\, A \neq 0
$$
and the corresponding irrationality measure function
$$
\Psi_\Lambda (t) = \min_{x= (x_1x_2) \in \Lambda: 0< \max(|x_1|,|x_2|)\le t}
\,\,\,\,
|x_1 x_2|^{1/2}.
$$
Following \cite{G2} we define the {\it ordinary Diophantine exponent of  a lattice}
$$
\omega_{\Lambda} = \sup\{\gamma\in \mathbb{R}: \liminf_{t\to \infty} t^\gamma \cdot 
\Psi_\Lambda (t)<\infty\}
$$
and the {\it weak uniform  Diophantine exponent of lattice}
$$
\overline{
\omega}_{\Lambda} =
 \sup\{\gamma\in \mathbb{R}: \limsup_{t\to \infty} t^\gamma \cdot 
\Psi_\Lambda (t)<\infty\}.
$$
Of course
\begin{equation}\label{tri}
\overline{
\omega}_{\Lambda} \le {
\omega}_{\Lambda}.
\end{equation}
We should notice that for any diagonal matrix
$$
D =
\left(
\begin{array}{cc}
d_1& 0\cr
0&d_1
\end{array}
\right),\,\,\,\,\,
d_1,d_2 \neq 0
$$
we have equalities
\begin{equation}\label{ekva}
\omega_{\Lambda_A} =\omega_{\Lambda_{DA}},\,\,\,\,\ 
\overline{\omega}_{\Lambda_A} =\overline{\omega}_{\Lambda_{DA}}.
\end{equation}
Let numbers 
\begin{equation}\label{the}
\theta = \frac{a_{1,2}}{a_{1,1}}, \,\,\, \eta = \frac{a_{2,1}}{a_{2,2}}
\end{equation}
be irrational.
 It is shown in \cite{G2} that  equalities (\ref{ekva}) immediately lead to
 \begin{equation}\label{matr}
 \omega_{\Lambda_A}
 =
 \frac{\max (\omega(\theta),\omega(\eta))+1}{2},\,\,\,\,\,
 \overline{\omega}_{\Lambda_A}=
  \frac{ {\varpi}_\upsilon (\theta,\eta) +1}{2}.
 \end{equation}
 We see that the study of  the connection between exponents $ \omega_{\Lambda_A},  \overline{\omega}_{\Lambda_A}$
 is reduced to the study of exponents $\omega(\theta),\omega(\eta)$ and  ${\varpi}_\upsilon (\theta,\eta)$.

\vskip+0.3cm
{\bf 2. Inequalities between Diophantine exponents}

\vskip+0.3cm
 Our first result is very easy.
  
  \vskip+0.3cm
{\bf Theorem 1.} {\it For any irrational $\theta$ one has
\begin{equation}\label{label}
\omega (\theta) 
\ge
\frac{1}{2-\overline{\omega}(\theta)},
\end{equation}
in particular, $\omega (\theta) =+\infty$ when $\overline{\omega}(\theta)=2$.
} 
\vskip+0.3cm

Proof. We may assume that $\overline{\omega}(\theta)>1$. Take
$ 1<\gamma< \overline{\omega}(\theta)$. Then for $ t$ large enough  by (\ref{alt}) we have
$$
\upsilon_\theta (t) = 
 \min_{\nu: \,\, q_\nu \le t } \,\, q_\nu\cdot ||q_\nu\theta||< t^{1-\gamma}.
$$
Let $\nu_*$ be defined by  inequalities $ q_{\nu_*}\le  t < q_{\nu_*+1}$.
Take $ t= q_{\nu_*}-1$ and let 
$$
\upsilon_\theta (q_{\nu_*}-1) =
 \min_{\nu<\nu_*} \,\, q_\nu\cdot ||q_\nu\theta|| =  q_{\nu_{**}}\cdot ||q_{\nu_{**}}\theta||,\,\,\,\,\,
 \nu_{**} <  \nu_{*}.
 $$
 By (\ref{2m}) we have
 $$
 a_{\nu_{**}+1} + 2 \ge (q_{\nu_*}-1)^{\gamma-1} \ge (q_{\nu_{**}+1}-1)^{\gamma-1}  \ge  (a_{\nu_{**}+1}q_{\nu_{**}})^{\gamma-1} ,\,\,\,\,\,
 \text{and so}\,\,\,\,\,
  a_{\nu_{**}+1}\gg q_{\nu_{**}}^{\frac{\gamma-1}{2-\gamma}}.
  $$
  Now
  $$
  || q_{\nu_{**}}\theta|| \le \frac{1}{q_{\nu_{**}+1}}\le \frac{1}{a_{\nu_{**}+1}q_{\nu_{**}}}\ll \frac{1}{q_{\nu_{**}}^{\frac{1}{2-\gamma}}},
  $$
  and everything is proven. $\Box$
 
 \vskip+0.3cm
 
 Next two results deal with  ordinary approximations to two real numbers.
 
 \vskip+0.3cm
 
{\bf Theorem 2.} {\it Assume that for irrational $\theta,\eta$ inequality $\varpi_\psi (\theta,\eta) >1$
holds.Then
$$
\min 
(
\omega (\theta),
\omega (\eta)
)
\ge  \varpi_\psi (\theta,\eta) ^2
$$

}

  \vskip+0.3cm
  To formulate a result dealing with weak approximations to two real numbers  we
  consider  polynomial
\begin{equation}\label{poly}
G_{y}(x) =
 x^2 -(y^2- 2 y+3) x +1.
\end{equation}
 It is clear  that 
$$
{g}(y) >1\,\,\,\,\,
\text{for}
\,\,\,\,\, y >1.
$$
We should notice that 
$$
G_{y}\left( y\right) = (1-y)^3\,\,\,\,\,
\text{and}
\,\,\,\,\,
G_{y}\left( \frac{1}{2-y}\right) = \left(\frac{y-1}{2-y}\right)^2.
$$
So for $y>1$  polynomial $G_y(x)$ has the lagrest root
$ \frak{g}_y , G_y(\frak{g}_y) = 0$ such that 
$$
\frak{g}_y >y.
$$
Moreover for $ 1<y<2$ this root satisfies the inequality
\begin{equation}\label{y}
\frak{g}_y <  \frac{1}{2-y}.
\end{equation}

\vskip+0.3cm

{\bf Theorem 3.} {\it Let
${\varpi}_\upsilon = {\varpi}_\upsilon (\theta,\eta)>1 $ and 
 $\frak{g}_{{\varpi}_\upsilon}$ be the  largest root of  polynomial $ G_{\varpi_\upsilon}(x)$.
Then 
\begin{equation}\label{ma1}
\max 
(
\omega (\theta),
\omega (\eta)
)
\ge \frak{g}_{{\varpi}_\upsilon}>  {\varpi}_\upsilon.
\end{equation}
%and

%НЕ ДОКАЗАНО В СЛУЧАЕ 1:

%\begin{equation}\label{ma2}
%\min 
%(
%\omega (\theta),
%\omega (\eta)
%)
%\ge {g}({{\varpi}_\upsilon})>  1.
%\end{equation}
}

\vskip+0.3cm
Define polynomial 
\begin{equation}\label{poly}
\mathscr{G}_y(x) = x^2 - 2yx-1.	 
\end{equation}
We see that for any $y>0$ there exists the root $\frak{G}_y$ of this polynomial  satisfying
$$
\mathscr{G}_y(x)=0,\,\,\,\,\, \frak{G}_y>y>1.
$$ 

  The following result about exponents of lattices
 is an immediate corollary of Theorem 3 and formulas (\ref{matr}).

\vskip+0.3cm
{\bf Theorem 4.} {\it Let $\Lambda= \Lambda_A$.
 Assume that numbers $\theta$ and $\eta$  defined in (\ref{the}) are irrational. Then
$$
\omega_\Lambda
\ge
\frak{G}_{\overline{\omega}_\Lambda}\cdot \overline{\omega}_\Lambda.
$$

}

 \vskip+0.3cm

 {\bf 4. Lemma about two piecewise constant functions.}
  \vskip+0.3cm
  The following lemma will be used in proofs of Theorems 2 and 3.

 \vskip+0.3cm
 
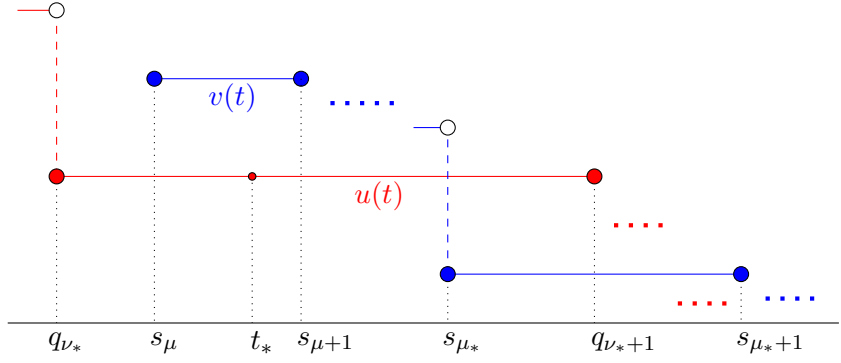
\begin{figure}[h]
  \centering
  \begin{tikzpicture}[scale=1.3]
    
    \draw[color=black] (-0.5,0) -- (8,0);
    
    \draw[color=red] (0,1.5) -- (5.5,1.5);
           \draw (3.3,1.3) node{\begin{small}{\color{red} $u (t)$}\end{small}};
            \draw (1.8,2.3) node{\begin{small}{\color{blue} $v (t)$}\end{small}};
        
    \draw[color=red] (-0.4,3.2) -- (0,3.2);

     \draw[dashed,color=red] (0,1.5) -- (0,3.2);
     
        \draw[dotted,color=black] (0,1.5) -- (0,0);
    
    \draw[color=blue] (1,2.5) -- (2.5,2.5);
        \draw[ultra thick,loosely dotted,color=blue] (2.8,2.25) -- (3.5,2.25);
        \draw[color=blue] (3.65,2) -- (4,2);
     \draw[dotted,color=black] (1,2.5) -- (1,0);

     \draw[color=blue] (7,0.5) -- (4,0.5);
      \draw[dotted,color=black] (4,0.5) -- (4,0);
        \draw[dashed,color=blue] (4,0.5) -- (4,2);

         \draw[dotted,color=black] (5.5,1.5) -- (5.5,0);
         
          \draw[dotted,color=black] (7,0.5) -- (7,0);
          
            \draw[dotted,color=black] (2,1.5) -- (2,0);
            
            \node[draw=black,fill=red,circle,inner sep=1pt] at (2,1.5) {};
            \node[draw=black,fill=red,circle,inner sep=2pt] at (0,1.5) {};
            \node[draw=black,fill=red,circle,inner sep=2pt] at (5.5,1.5) {};
            
                \node[draw=black,fill=blue,circle,inner sep=2pt] at (2.5,2.5) {};
                
                  \draw[dotted,color=black] (2.5,2.5) -- (2.5,0);
                   \node[draw=black,fill=blue,circle,inner sep=2pt] at (1,2.5) {};
                   \node[draw=black,fill=white,circle,inner sep=2pt] at (4,2) {};
                   \node[draw=black,fill=blue,circle,inner sep=2pt] at (4,0.5) {};
                   \node[draw=black,fill=blue,circle,inner sep=2pt] at (7,0.5) {};
                   
                   \node[draw=black,fill=white,circle,inner sep=2pt] at (0,3.2) {};

                   \draw[ultra thick,loosely dotted,color=blue] (7.25,0.25) -- (7.75,0.25);
                   
                    \draw[ultra thick,loosely dotted,color=red] (5.7,1) -- (6.3,1);
                    
                        \draw[ultra thick,loosely dotted,color=red] (6.85,0.2) -- (6.35,0.2);

                      \draw (0.1,-0.2) node{\begin{small}$q_{\nu_*}$\end{small}};
                      
                      \draw (1.1,-0.2) node{\begin{small}$s_{\mu}$\end{small}};
                      
                       \draw (2.1,-0.2) node{\begin{small}$t_*$\end{small}};
                       
                        \draw (2.75,-0.2) node{\begin{small}$s_{\mu+1}$\end{small}};
                        
                                     \draw (4.15,-0.2) node{\begin{small}$s_{\mu_*}$\end{small}};

                                     \draw (7.3,-0.2) node{\begin{small}$s_{\mu_*+1}$\end{small}};
                                     
                                      \draw (5.8,-0.2) node{\begin{small}$q_{\nu_*+1}$\end{small}};
    
  \end{tikzpicture}
  
     \caption{to the proof  of Lemma 1} \label{fig0}

  \end{figure}

 {\bf Lemma 1.} {\it
 Consider two  increasing sequences
 $$
 q_1<q_2<...<q_\nu<q_{\nu+1}<...\, ,\,\,\,\,\,\,\,
 s_1<s_2<...<s_\mu<s_{\mu+1}<...
 $$ and two positive-valued piecewise constant decreasing functions
 $$
 u(t) = u(q_\nu)\,\,\,\text{for} \,\,\, t\in [q_\nu,q_{\nu+1});\,\,\,\,\,  u(q_\nu)\ge  u(q_{\nu+1})\,\,\,\text{for every}\,\,\,\nu
 $$
 and
  $$
 v(t) = v(s_\mu)\,\,\,\text{for} \,\,\, t\in [s_\mu,s_{\mu+1});\,\,\,\,\,  v(s_\mu)\ge  v(s_{\mu+1}) \,\,\,\text{for every}\,\,\,\mu.
 $$
 Assume that 

 \noindent
{\rm ({\bf a})  }
for every $\mu$  large enough there exists $ t_{*} \in [s_\mu, s_{\mu+1})$ such that  the strict inequality
$$
 u  (t_{*}) < v (t_{*})
$$
holds
and

 \noindent
{\rm ({\bf b})  }
for every $\nu $  large enough there exists $ t_{**} \in [q_\nu, q_{\nu+1})$ such that   the strict inequality
$$
 v  (t_{**}) < u (t_{**})
$$
holds.
Then there exist arbitrary large values of indices $\nu_*$ and $\mu_*$ such that 
\begin{equation}\label{e00}
 u  (t)
 \,\,\,
 \text{is not continuous at point}
 \,\,\,
 t = q_{\nu_*},
 \end{equation}
 \begin{equation}\label{e0000}
 \upsilon  (t)
 \,\,\,
 \text{is not continuous at point}
 \,\,\,
 t = s_{\mu_*},
\end{equation} 
 \begin{equation}\label{ee0000}
q_{\nu_*}<s_{\mu_*}<q_{\nu_*+1}< s_{\mu_*+1},
\end{equation}
and

\begin{equation}\label{1w}
 u(s_{\mu_*})< v(s_{\mu_*}-),
\end{equation}
\begin{equation}\label{2w}
 v(q_{\nu_*+1}-)< u(q_{\nu_*+1}-)
\end{equation}

}
 \vskip+0.3cm
The proof of Lemma 1 is visualised at Fig. \ref{fig0}.

\vskip+0.3cm
 
 Proof.
 Fix $\mu $ and
 consider $\nu_*$ such that for $t_*$ from  condition ({\bf a})   one has $ t_* \in [q_{\nu_*}, q_{\nu_*+1})$. Then obviously
 \begin{equation}\label{0e}
 q_{\nu_*}\le t_* < s_{\mu+1}
 \end{equation}
 and
\begin{equation}\label{0ee}
u (q_{\nu_*})= u (t_*) <v (t_{*}).
\end{equation}
Moreover, from condition ({\bf b}) it follows that  
\begin{equation}\label{0eee}
u(q_{\nu_*-1}) = u  (q_{\nu_*})
\end{equation}
 is not possible. Indeed in the case of equality 
 (\ref{0eee}) 
 by (\ref{0ee}) 
 for every $  t\in [q_{\nu_*-1}, q_{\nu_*}) \subset [ q_{\nu_*-1}, t_*]$ 
  we get
 $$
 u(t) =
 u  (q_{\nu_*})
<  v (t_{*}) \le  v  (t) 
 $$
 what is a contradiction to condition ({\bf b}).
 So  we have a strict inequality  $u(q_{\nu_*-1}) >  u(q_{\nu_*})$
and (\ref{e00}) is valid.

Then, from condition ({\bf b})  
for $\nu = \nu_*$ we see that 
$$
 s_{\mu+1} <q_{\nu_*+1}. 
$$
Now we take $\mu_*\ge\mu+1$ such that $s_{\mu_*}\le q_{\nu_*+1}< s_{\mu_*+1}$. 
We should notice that 
\begin{equation}\label{ee00e}
s_{\mu+1}\le s_{\mu_*}.
\end{equation}
We show that the first inequality here is a strict one, that is
\begin{equation}\label{ee00}
s_{\mu_*}<q_{\nu_*+1}< s_{\mu_*+1}.
\end{equation}
Indeed, in the case of equality $s_{\mu_*}=q_{\nu_*+1}$ there may be two opportunities.

The {\bf first} opportunity is 
$$ v(s_{\mu_*}-)\ge u(s_{\mu_*}-) = u(q_{\nu_*+1}-). $$
But in such a case for any $t \in [q_{\nu_*}, q_{\nu_*+1}) =  [q_{\nu_*}, s_{\mu_*})$ we have 
$$
u(t) =  u(q_{\nu_*+1}-)\le v(s_{\mu_*}-)\le v(t),
$$
and this contradicts condition ({\bf b}) for $ \nu = \nu_*$.

The {\bf second} opportunity is 
$$ v(s_{\mu_*}-) < u(s_{\mu_*}-) = u(q_{\nu_*+1}-). $$
We consider $s_{\mu_*-1}< s_{\mu_*}  =q_{\nu_*+1}$ and conclude 
  that  for any $t \in [  s_{\mu_*-1}, q_{\nu_*+1})$ one has
\begin{equation}\label{ee00ee}
v(t) =
v(s_{\mu_*-1}) =v(s_{\mu_*}-) <  u(q_{\nu_*+1}-) \le u(t). 
\end{equation}
This is a contradiction to condition ({\bf a}) with $\mu = \mu_*-1$. So we have proven (\ref{ee00}).

Now we see that inequalities (\ref{0e},\ref{ee00e},\ref{ee00}) lead to (\ref{ee0000}).

Next we prove (\ref{e0000}).   For the value $v(s_{\mu_*})$ we again should consider two opportunities.

The {\bf first} opportunity is  
\begin{equation}\label{opa}
 v(s_{\mu_*})= 
 v(s_{\mu_*+1}-)
\ge u(q_{\nu_*+1}-).
\end{equation}
As 
$$
v(t) \ge   v(s_{\mu_*}) \,\,\,\,\text{for all}\,\,\,t < s_{\mu_*+1};\,\,\,\,\,
u(t) = u(q_{\nu_*+1}-) \,\,\,\,\text{for all}\,\,\,t \in [q_{\nu_*},q_{\nu_*+1})
$$
and
$
q_{\nu_*+1}<s_{\mu_*+1}
$
we immediately get a contradiction with condition ({\bf b}) for $\nu = \nu_*$. So the  {\bf first} opportunity is not possible.

The {\bf second} opportunity is 
$$ v(s_{\mu_*}) 
< u(q_{\nu_*+1}-).
$$
But if $v(t)$ is continuous at  point $ t = s_{\mu_*}$ then $v(s_{\mu_*-1})= v(s_{\mu_*})$ and we again have inequality
(\ref{ee00ee}) which leads to a contradiction.
So (\ref{e0000}) is proven.

Now we prove (\ref{1w}). Assuming that (\ref{1w}) is not valid,   for every $ t\in [s_{\mu_*-1},s_{\mu_*}) $ we immediately obtain inequality
$$
u(t) \ge u(s_{\mu_*}) \ge  v(s_{\mu_*}-)  =   v(s_{\mu_*-1}) = v(t)
$$
which contradicts to condition ({\bf a}) with $ \mu = \mu_*-1$.
So (\ref{1w}) is proven

Finally,
as (\ref{opa}) is not possible,
we have 
$$
v(q_{\nu_*+1}-) \le v(s_{\mu_*}) < u(q_{\nu_*+1}-),
$$
and this proves (\ref{2w}).$\Box$

 \vskip+0.3cm

{\bf 4. Proofs of Theorems 2,3.}

 \vskip+0.3cm
 Proofs of Theorems 2 and  3  rely on Lemma 1.  We apply this lemma to
 functions
  \begin{equation}\label{fff0}
 v(t) = \psi_\eta(t),\,\,\,\,\,\,\,\,\,\,
 u(t) =\psi_\theta(t)
 \end{equation}
 to get a proof of Theorem 2 and to functions 
 \begin{equation}\label{fff}
 v(t) = \upsilon_\eta(t),\,\,\,\,\,\,\,\,\,\,
 u(t) =\upsilon_\theta(t)
 \end{equation}
 to  get a proof of Theorem 3.
 As for the sequences involved in Lemma 1, we  consider sequence $\frak{Q}$ o
 of convergents' denominators for  $\theta$
  introduced in (\ref{1mx})  in the beginning of the paper and sequence 
 $$
\frak{S}:
\,\,\,\,\,
s_1<s_2<...<s_\mu<s_{\mu+1}<...
$$
of convergents' denominators for  $\eta$.

  \vskip+0.3cm

  \vskip+0.3cm
 {\bf 3.1. Proof of Theorem 2.}
 
   \vskip+0.3cm

We  want to apply  Lemma 1  to functions  (\ref{fff0}) and sequences $\frak{Q},\frak{S}$.
First of all we should show that  conditions ({\bf a}) and ({\bf b}) are satisfied under the assumptions of Theorem 2.
Take $\gamma \in (1,\varpi ({\theta,\eta}))$.   Then by definition of $\varpi_\psi (\theta,\eta)$,
for large values of $t$ we have
$$
\min ( \psi_\theta (t) , \psi_\eta (t)) < \frac{1}{t^\gamma}.
$$
But 
for any values of $\nu, \mu$ one has
$$
q_{\nu+1}\psi_\theta (q_{\nu+1}-) =
q_{\nu+1}||q_\nu \theta|| >\frac{1}{2},\,\,\,\,\,
s_{\mu+1}\psi_\eta (s_{\nu+1}-) =
s_{\mu+1}||s_\mu \eta|| >\frac{1}{2}.
$$ 
If   ({\bf a}) or  ({\bf b}) is not satisfied, then either there exist arbitrary large $\nu$ with 
$$
\min ( \psi_\theta (t) , \psi_\eta (t)) = \psi_\theta (t) \,\,\,\,\,
\text{for all}\,\,\, t \,\,\,
\text{from the interval}
\,\,\,\,\, q_\nu \le t < q_{\nu+1}
,
$$
or there exist arbitrary large $\mu$ with 
$$
\min ( \psi_\theta (t) , \psi_\eta (t)) = \psi_\eta (t) \,\,\,\,\,
\text{for all}\,\,\, t \,\,\,
\text{from the interval}
\,\,\,\,\, s_\mu \le t < s_{\mu+1}.
$$
In both cases we have a contradiction.
So we may apply Lemma 1.

Let $\nu_*,\mu_*$ be as in Lemma 1.
The very
first inequality from (\ref{2m}) applied to approximations to $\theta$, formula (\ref{1w})  and definition of $\varpi ({\theta,\eta})$ show that for large values of indices one has
$$
\frac{1}{2q_{\nu_*+1}} \le \psi_\theta (q_{\nu_*}) = \psi_\theta (s_{\mu_*}-) =
\psi_{\theta,\eta} (s_{\mu_*}-) \le \frac{1}{s_{\mu_*}^\gamma}.
$$
Similarly
first inequality from (\ref{2m}) applied to approximations to $\eta$, formula (\ref{2w})  and definition of $\varpi_\psi ({\theta,\eta})$ show that for large values of indices one has
$$
\frac{1}{2s_{\mu_*+1}} \le  \psi_\eta (s_{\mu_*})=
\psi_\eta (q_{\nu_*+1}-) = \psi_{\theta, \eta} (q_{\nu_*+1}-) \le \frac{1}{q_{\nu_*+1}^\gamma}.
$$
 From the last two inequalities we see that 
 $$
s_{\mu_*+1}\ge \frac{s_{\mu_*}^{\gamma^2}}{2^{1+\gamma}}
$$
and by the second inequality from   (\ref{2m})  we have
$$
||s_{\mu_*}\eta||\le \frac{1}{s_{\mu_*+1}}\ll \frac{1}{s_{\mu_*}^{\gamma^2}}.
$$
So we proved $\omega (\eta) \ge\varpi^2_\psi ({\theta,\eta})$. To prove  $\omega (\theta) \ge\varpi^2_\psi ({\theta,\eta})$ we need to exchange the roles of
$\theta$ and $\eta $ in our proof. 
$\Box$.
 
\vskip+0.3cm

{\bf 3.2. Proof of Theorem 3.}
\vskip+0.3cm

Take $\gamma < {\varpi}_\upsilon-1$ close to ${\varpi}_\upsilon-1$.
We consider two cases.

\vskip+0.3cm
{\bf Case 1$^0$}.   Either

\noindent
({\bf 1a})
 There exists infinitely many  consecutive elements  $q_\nu,q_{\nu+1}$ of $\frak{Q}$ such that 
$$
 \upsilon_{\theta,\eta}  (t)=   \upsilon_{\theta}  (t)\,\,\,\,\,\text{for all}\,\,\, t\,\,\,
 \text{from the interval}\,\,\,
[q_\nu,q_{\nu+1}),
$$
 or
 
 \noindent
({\bf 1b})
  similarly there exists infinitely many  consecutive elements  $
s_\mu,s_{\mu+1}$ of $\frak{S}$ such that 
$$
 \upsilon_{\theta,\eta}  (t)=   \upsilon_{\eta}  (t)\,\,\,\,\,\text{for all}\,\,\, t\,\,\,
 \text{from the interval}\,\,\,
[s_\mu,s_{\mu+1}).
$$

\vskip+0.3cm
Assume that ({\bf 1a}) holds (under assumption  ({\bf 1b})  the argument is quite similar). Then by  lower bound from (\ref{2m})  and  (\ref{alt}) we get
$$
\min_{j\le \nu} \,\, \frac{1}{a_{j+1}+2} \le 
\min_{j\le \nu} \,\, q_j ||q_j\theta|| = \upsilon_\theta (q_{\nu+1}-) =
\upsilon_{\theta,\eta} (q_{\nu+1}-) 
 \le q_{\nu+1}^{-\gamma}.
$$
As $ a_{j+1}+2< q_{\nu+1}$ for large  $j\le \nu$ we see that $\gamma <1$. This means that  for $ {\varpi}_\upsilon> 2$ and
$ \gamma$ close to $ {\varpi}_\upsilon$ {\bf case 1$^0$}  cannot happen.
 
 Now consider $j_*$ such that 
 $$
a_{j_*+1}=
\max_{j\le \nu }\,\, a_{j+1} \gg q_{\nu+1}^\gamma  \ge q_{j_*+1}^\gamma \gg
(a_{j_*+1}q_{j_*})^\gamma.
$$
Then
$$
a_{j_*+1}\gg q_{j_*}^{\frac{\gamma}{1-\gamma}}
.
$$
Finally,
$$
||q_{j_*}\theta|| \le \frac{1}{a_{j_*+1}q_{j_*}}\ll  q_{j_*}^{-\frac{1}{1-\gamma}}.
$$
In such a way by (\ref{y})  in {\bf case 1$^0$}   under condition
({\bf 1a})
we get $\omega(\theta) \ge \frac{1}{2-{\varpi}_\upsilon} \ge \frak{g}_{{\varpi}_\upsilon}$. 
Of course, under condition
({\bf 1b})
we get $\omega(\eta) \ge \frac{1}{2-{\varpi}_\upsilon} \ge \frak{g}_{{\varpi}_\upsilon}$.

\vskip+0.3cm
We repeat that {\bf case 1$^0$} deals with the situation when there are infinitely many consecutive convergents' denominators for $\theta$ or $\eta$ such that in the whole interval between these two consecutive denominators  the minimum 
$
\upsilon_{\theta,\eta}  (t)= \min ( \upsilon_\theta (t), \upsilon_\eta (t))
$
is equal to one of the functions  
$ \upsilon_\theta (t)$ or $ \upsilon_\eta (t)$ (namely just  to the function
corresponding to the number with the consecutive denominators under the consideration).
The alternative to  {\bf case 1$^0$} is the following
\vskip+0.3cm
{\bf Case 2$^0$}. There exist  $\nu_0,\mu_0$ such that 
simultaneously

 \noindent
({\bf 2a})  
for every $\mu \ge \mu_0$ there exists $ t_{*} \in [s_\mu, s_{\mu+1})$ such that 
$$
 \upsilon_{\theta}  (t_{*}) <  \upsilon_{\eta}  (t_{*})
$$
and

 \noindent
({\bf 2b})  
for every $\nu \ge \nu_0$ there exists $ t_{**} \in [q_\nu, q_{\nu+1})$ such that 
$$
 \upsilon_{\eta}  (t_{**}) <  \upsilon_{\theta}  (t_{**}).
$$
\vskip+0.3cm

 In {\bf case 2$^0$} 
 conditions ({\bf a}) and ({\bf b})  
 of Lemma 1
 for functions (\ref{fff})
  are satisfied.

\vskip+0.3cm
Now we combine the last  inequality from (\ref{2m}) in its asymptotic form with (\ref{e0000}) and  (\ref{alto}) applied to  the function 
 $\upsilon_{\eta}  (t)$ at point   $ t = s_{\mu_*}$, and  take into account  (\ref{2w}) to get 
\begin{equation}\label{X0}
\frac{s_{\mu_*}}{s_{\mu_*+1}}\asymp  \upsilon_\eta (s_{\mu_*}) =
\upsilon_\eta (q_{\nu_*+1}- )=
\upsilon_{\theta,\eta} (q_{\nu_*+1}- )
 \le \frac{1}{q_{\nu_*+1}^\gamma}
\end{equation}
Similarly, combining
  inequality (\ref{2m}) in its asymptotic form with (\ref{e00}) and  (\ref{alto}) applied to  the function 
 $\upsilon_{\theta}  (t)$ at point   $t = q_{\nu_*}$, and taking into account  (\ref{1w}), we get 
\begin{equation}\label{X00}
\frac{q_{\nu_*}}{q_{\nu_*+1}} \asymp \upsilon_\theta (q_{\nu_*}) 
= \upsilon_\theta (s_{\mu_*}) =  \upsilon_{\theta,\eta} (s_{\mu_*}-) 
\le \frac{1}{s_{\mu_*}^\gamma}.
\end{equation}
Finally, inequalities (\ref{X0}) and (\ref{X00}) can be rewritten as
\begin{equation}\label{X000}
q_{\nu_*+1}^\gamma s_{\mu_*} \ll s_{\mu_*+1},\,\,\,\,\,
q_{\nu_*} s_{\mu_*}^\gamma \ll q_{\nu_*+1}.
\end{equation}
\vskip+0.3cm
Now we are able to finalise the proof of Theorem 3.

\vskip+0.3cm
First of all, assume that 
$$
q_{\nu_*+1}\le q_{\nu_*}^g,\,\,\,\,\,
s_{\mu_*+1}\le s_{\mu_*}^g\,\,\,\,
\,
\text{with some }
\,\,\,\,\, g>1.
$$
Then the first inequality from (\ref{X000}) gives
$$
q_{\nu_*+1}\ll s_{\mu_*}^{\frac{g-1}{\gamma}},
$$
meanwhile the second inequality from (\ref{X000}) gives
$$
s_{\mu_*}\ll q_{\nu_*+1}^{\frac{1}{\gamma} \left( 1- \frac{1}{g}\right)}.
$$
Combining the last two inequalities we see that 
$$
(g-1)^2\ge \gamma^2g,
$$
and this means that $ g$ cannot be smaller than the root $\frak{g}_{\gamma+1}$ of the polynomial
$ G_{
\gamma+1}(x)
= (g-1)^2 -\gamma^2 g
$ 
defined in (\ref{poly}).  So
$$
\text{either}\,\,\,\,\,
q_{\nu_*+1}\ge q_{\nu_*}^{\frak{g}_{\gamma+1}},\,\,\,\,\,
\text{or}\,\,\,\,\,
s_{\mu_*+1}\ge s_{\mu_*}^{\frak{g}_{\gamma+1}}.
$$
By the upper bound from (\ref{2m}),  we see that 
$$
\text{either}\,\,\,\,\,
||q_{\nu_*}\theta||\le q_{\nu_*}^{-\frak{g}_{\gamma+1}},\,\,\,\,\,
\text{or}\,\,\,\,\,
||s_{\mu_*}\eta||\le s_{\mu_*}^{-\frak{g}_{\gamma+1}}.
$$
As we can take $\gamma$ arbitrary close to $\overline{\varpi}-1$, this means that 
$$
\text{either}\,\,\,\,\,
 \omega(\theta) \ge 
 \frak{g}_{{\varpi}_\upsilon},\,\,\,\,\,
\text{or}\,\,\,\,\,
 \omega(\eta) \ge 
 \frak{g}_{{\varpi}_\upsilon},
$$
and (\ref{ma1}) is proven.$\Box$
\vskip+0.3cm

%As for (\ref{ma2}), we will prove here  that 
%\begin{equation}\label{y000}
%\omega (\theta) \ge g({\varpi}_\upsilon).
%\end{equation}
 %To prove inequality 
%$ \omega (\eta) \ge g(\overline{\varpi})$ one should exchange the roles of $\theta$ and $\eta$ in the whole  proof from %the very beginning.

%Indeed, assume that 
%$$
%s_{\mu_*+1} \le s_{\mu_*}^g\,\,\,\,\text{with}\,\,\,\, g>1.
%$$
%From  both inequalities  (\ref{X000}) we deduce
%$$
%s_{\mu_*}^{\gamma^2+1} = (s_{\mu_*}^\gamma)^\gamma\le
% (q_{\nu_*}s_{\mu_*}^\gamma)^\gamma
 %\ll q_{\nu_*+1}^\gamma s_{\mu_*}\ll  s_{\mu_*+1}\le  s_{\mu_*}^g,
 %$$
 %and so $ g\ge \gamma^2+1$. This is true for any $ \gamma < \overline{\varpi}-1$. The last inequality gives (\ref{y000}).$\Box$

 \vskip+0.3cm
 
{\bf 4. About optimality of inequalities} 
\vskip+0.3cm
By constructing examples of $\theta,\eta$ we will show here that the inequalities of Theorems 1, 2, 3 are optimal.

\vskip+0.3cm

{\bf 4.1.} Optimality of  inequality (\ref{label}) of Theorem 1 is almost obvious. For $\gamma \in (1,2)$ we take $\theta$ with partial quotients defined inductively by $ a_{\nu+1} \asymp q_\nu^{\frac{\gamma -1}{2-\gamma}}$. Then
$$
q_{\nu+1}\asymp q_\nu^\frac{1}{2-\gamma},\,\,\,\,\,
\upsilon_\theta (q_\nu) = \upsilon_\theta(q_{\nu+1}-) = q_\nu ||q_\nu \theta||\asymp \frac{1}{a_{\nu+1}} \asymp \frac{1}{q_{\nu+1}^{\gamma -1}}
$$
and so $ \omega(\theta) = \frac{1}{2-\gamma}$ and $ \overline{\omega} (\theta) = \gamma$.

\vskip+0.3cm

{\bf 4.2.}  To show optimality of Theorem 2  we take $\gamma >1$ and construct  numbers $\theta$ and $\eta $
by defining partial quotients $a_\nu, b_\mu$ in their continued fractions' expansions
\begin{equation}\label{rt1}
\theta=
   [a_0;a_1,a_2,...,a_\nu,...],\,\,\,\,\,
   \eta=
   [b_0;b_1,b_2,...,b_\mu,...].
\end{equation}
We can inductively choose $a_\nu, b_\mu$ to satisfy
\begin{equation}\label{rt2}
q_{\nu}  \sim s_{\nu}^\gamma,\,\,\,\,\,
s_{\nu+1} = b_{\nu+1}s_\nu + s_{\nu-1} \sim q_{\nu}^\gamma,\,\,\,\,\,
q_{\nu+1}  = a_{\nu+1}q_\nu + a_{\nu-1} \sim s_{\nu+1}^\gamma,\,\,\,\,\,
\nu \to \infty
.
\end{equation}
Then
\begin{equation}\label{rt23}
s_1<q_1<...<s_\nu<q_\nu< s_{\nu+1}<q_{\nu+1}<...
\end{equation}
and
$$
q_{\nu+1} \sim q_\nu^{\gamma^2},\,\,\,\,\,\,\,\,
s_{\nu+1} \sim s_\nu^{\gamma^2}
.
$$
We see that 
$$
\psi_\theta(s_{\nu}) =
\psi_\theta(q_{\nu-1}) \asymp \frac{1}{q_{\nu}}
<
\psi_\eta(s_{\nu}) \asymp \frac{1}{s_{\nu+1}} \asymp \frac{1}{q_{\nu}^\gamma}
$$
and similarly
$$
\psi_\theta(q_\nu) \asymp \frac{1}{q_{\nu+1}}  \asymp \frac{1}{s_{\nu+1}^\gamma} 
< \psi_\eta(q_\nu) = \psi_\eta(s_\nu)\asymp \frac{1}{s_{\nu+1}}.  
$$
So
$$
\psi_{\theta,\eta}(t) =
\min
(
\psi_{\theta}(t) , 
\psi_{\eta}(t) )=
\begin{cases}
\psi_{\eta}(s_\nu) \asymp\frac{1}{s_{\nu+1}}\asymp \frac{1}{q_\nu^\gamma},\,\,\,\,\,\,\,\,\text{if}\,\,\,\,\,
t \in [s_\nu, q_\nu),
\cr
\psi_{\theta}(q_\nu) \asymp\frac{1}{q_{\nu+1}}\asymp \frac{1}{s_{\nu+1}^\gamma},\,\,\,\,\text{if}\,\,\,\,\,
t \in [q_\nu, s_{\nu+1}).
\end{cases}
$$
Finally  
\begin{equation}\label{rt}
\psi_{\theta,\eta}(t) \ll t^{-\gamma}\,\,\,\,\,
\text{for all}\,\,\,\,\, t
.
\end{equation}
So,  by  construction we have $ \omega(\theta) =\omega (\eta) = \gamma^2$ and from Theorem 2 and (\ref{rt}) we conclude that 
$  \varpi_\psi (\theta,\eta) = \gamma$.

\vskip+0.3cm
{\bf 4.3.}  To show optimality of Theorem 3  we take $\gamma >0$ and construct  numbers $\theta$ and $\eta $ in the form 
(\ref{rt1}) with partial quotients  $a_\nu, b_\nu$ that ensure conditions
\begin{equation}\label{rt4}
s_{\nu+1} = q_{\nu}^\gamma s_\nu+ O(s_\nu),\,\,\,\, q_{\nu+1} =  q_\nu s_{\nu+1}^\gamma+O(q_\nu),\,\,\,\,\,
\nu \to \infty
\end{equation}
instead of (\ref{rt2}), and get sequences (\ref{rt23}). Asymptotic equalities (\ref{rt4}) show that

$$
\left(
\begin{array}{c}
\log q_{\nu+1}\cr
\log s_{\nu+1}
\end{array}
\right)
=
\left(
\begin{array}{cc}
\gamma^2+1& \gamma \cr
 \gamma & 1
\end{array}
\right)\left(
\begin{array}{c}
\log q_{\nu}\cr
\log s_{\nu}
\end{array}
\right)
+O\left(
s_\nu^{-\gamma}
\right)
.
$$
As the characteristic polynomial of the matrix here is
$$
x^2 - (\gamma^2+2)x +1 =(x-1)^2 - \gamma^2 x,
$$ 
we conclude that 
$$
\log q_\nu = A g^\nu+ O(1),\,\,\,\,\, 
\log s_\nu = B g^\nu+ O(1),
$$
where $ g>1$ is the root of the characteristic  polynomial and $A,B$ are positive constants. So
$$
 q_{\nu+1}\asymp q_\nu^g,\,\,\,\,\,\,\, s_{\nu+1}\asymp s_\nu^g
$$
and
$ \omega(\theta) =\omega (\eta) = g$.
At the same time
$$
\upsilon_{\theta,\eta}(t) =
\min
(
\upsilon_{\theta}(t) , 
\upsilon_{\eta}(t) )=
\begin{cases}
\upsilon_{\eta}(s_\nu) \asymp\frac{s_\nu}{s_{\nu+1}}\asymp \frac{1}{q_\nu^\gamma},\,\,\,\,\,\,\,\,\text{if}\,\,\,\,\,
t \in [s_\nu, q_\nu),
\cr
\upsilon_{\theta}(q_\nu) \asymp\frac{q_\nu}{q_{\nu+1}}\asymp \frac{1}{s_{\nu+1}^\gamma},\,\,\,\,\text{if}\,\,\,\,\,
t \in [q_\nu, s_{\nu+1}).
\end{cases}
$$
and this means that 
${\varpi}_\upsilon (\theta,\eta) = \gamma +1$.

\end{document}